\documentclass[12pt]{amsart}
\usepackage{psfig}

\newcommand{\R}{\mathbb{R}}
\newcommand{\Z}{\mathbb{Z}}

\newcommand{\N}{\mathbb{N}}
\newcommand{\A}{\mathbb{A}}
\newcommand{\T}{\mathbb{T}}

\newcommand{\F}{\mathcal{F}}
\newcommand{\G}{\mathcal{G}}

\newcommand{\M}{\mathcal{M}}
\newcommand{\g}{\gamma}
\newcommand{\e}{\epsilon}

\def\SBIMSMark#1#2#3{
 \font\SBF=cmss10 at 10 true pt
 \font\SBI=cmssi10 at 10 true pt
 \setbox0=\hbox{\SBF Stony Brook IMS Preprint \##1}
 \setbox2=\hbox to \wd0{\hfil \SBI #2}
 \setbox4=\hbox to \wd0{\hfil \SBI #3}
 \setbox6=\hbox to \wd0{\hss
             \vbox{\hsize=\wd0 \parskip=0pt \baselineskip=10 true pt
                   \copy0 \break%
                   \copy2 \break%
                   \copy4 \break}}
 \dimen0=\ht6   \advance\dimen0 by \vsize \advance\dimen0 by 8 true pt
                \advance\dimen0 by -\pagetotal
 \dimen2=\hsize \advance\dimen2 by .25 true in
  \openin2=publishd.tex
  \ifeof2\setbox0=\hbox to 0pt{}
  \else 
     \setbox0=\hbox to 3.1 true in{
                \vbox to \ht6{\hsize=3 true in \parskip=0pt  \noindent  
                \input publishd.tex 
                \vfill}}
  \fi
  \closein2
  \ht0=0pt \dp0=0pt
 \ht6=0pt \dp6=0pt
 \setbox8=\vbox to \dimen0{\vfill \hbox to \dimen2{\copy0 \hss \copy6}}
 \ht8=0pt \dp8=0pt \wd8=0pt
 \copy8
 \message{*** Stony Brook IMS Preprint #1, #2 ***}
}

\begin{document}
\SBIMSMark{1997/3}{February 1997}{}

\theoremstyle{plain}
\newtheorem{theorem}{Theorem}
\newtheorem{lemma}{Lemma}
\newtheorem*{condition}{}
\newtheorem{corollary}{Corollary}
\newtheorem*{thm1p}{Theorem 1${}^\prime$}
\theoremstyle{definition}
\newtheorem*{definition}{Definition}
\newtheorem*{remark}{Remark}
\theoremstyle{remark}

\title{Stably Non-Synchronizable Maps of the Plane}

\author[P.~Le~Calvez]{Patrice~Le~Calvez}
\address{URA 742 du CNRS, Institut Galil{\'e}e, Universit{\'e} Paris
13, Avenue J.-B. Cl{\'e}ment, 93420 Villetaneuse, France.}
\email{lecalvez@math.univ-paris13.fr}

\author[M.~Martens]{Marco~Martens}
\address{Institute of Mathematical Sciences, SUNY at
Stony Brook, Stony Brook, NY 11794-3651.}
\email{marco@math.sunysb.edu}

\author[C.~Tresser]{Charles~Tresser}
\address{I.B.M., P.O. Box 218, Yorktown Heights, NY 10598.}
\email{tresser@watson.ibm.com}

\author[P.~Worfolk]{Patrick~A.~Worfolk}
\address{The Geometry Center, 1300 S. Second St., Minneapolis,
MN 55454.}
\email{worfolk@geom.umn.edu}

\begin{abstract}
Pecora and Carroll presented a notion of synchronization where an
$(n-1)$-dimensional nonautonomous system is constructed from a given
$n$-dimensional dynamical system by imposing the evolution of one
coordinate.  They noticed that the resulting dynamics may be
contracting even if the original dynamics are not.  It is easy to
construct flows or maps such that no coordinate has synchronizing
properties, but this cannot be done in an open set of linear maps or
flows in $\R ^n$, $n\geq 2$.  In this paper we give examples of real
analytic homeomorphisms of $\R ^ 2$ such that the non-synchronizability
is stable in the sense that in a full $C^0$ neighborhood of the given
map, no homeomorphism is synchronizable.
\end{abstract}

\maketitle

%
%

\section{Introduction}

The observation that nonlinear oscillators can synchronize when
coupled was first reported in a letter that Christiaan Huyghens wrote
to his father \cite{Hu}. This phenomenon has inspired many
investigations by physicists and mathematicians over the
years. Another (although related) kind of synchronization has recently
received much attention because of its potential applications, and
because of the surprise caused by its discovery by Pecora and Carroll
\cite{PC1}.  Roughly speaking, they noticed that if one constructs an
$(n-1)$-dimensional nonautonomous evolution equation from a given
$n$-dimensional one by imposing the evolution of one coordinate, the
resulting dynamics could be contracting, even when the original
dynamics are not.  The surprise was that this works for chaotic
evolution equations such as the Lorenz equations as observed in
\cite{PC1} and subsequently proved in \cite{HV}. We first recall what
we need from \cite{PC1} and \cite{PC2}.

Consider the Lorenz equations:
$$
\left\{ \begin{array}{rl}
		\dot{x} &= \sigma (y-x)\,, \\
		\dot{y} &= rx-y-xz\,, \\
		\dot{z} &= xy-bz\,.
	 \end{array} \right.
$$
Pecora and Carroll \cite{PC1} noticed that for any solution
$(x(t),y(t),z(t))$ of these equations, all solutions $(Y(t),Z(t))$ of
the nonautonomous system
$$
\left\{ \begin{array}{rl}
	\dot{Y} &= rx(t)-Y-x(t)Z \,, \\
	\dot{Z} &= x(t)Y-bZ \,,
	\end{array}\right.
$$
satisfy
$$
\lim_{t\to+\infty}|y(t)-Y(t)|=\lim_{t\to+\infty}|z(t)-Z(t)| = 0 \,.
$$
They called this phenomenon {\it synchronization}, and the name {\it
master-slave synchronization} was proposed in \cite{TWB} to avoid
confusion with other, previously recognized, phenomena.  More is true
for the Lorenz system: for any smooth function $x(t)$, not necessarily
a solution of the Lorenz equations,
$$
	\lim_{t\to +\infty}|Y_1(t)-Y_2(t)| =
	\lim_{t\to +\infty}|Z_1(t)-Z_2(t)| = 0 \,,
$$
for any two initial conditions $(Y_1(0),Z_1(0))$, $(Y_2(0),Z_2(0))$ of
the above nonautonomous system. This stronger property was called {\it
absolute (master-slave) synchronization} in \cite{TWB} where all these
concepts were defined in a more geometric and coordinate independent
way.  Since the only kind of synchronization we will be dealing with
in this paper is of the master-slave kind, we shall from now on omit
this qualification.  Synchronization also occurs for maps, and this
will be our sole concern in this paper. Accordingly we next recall the
relevant definitions from \cite{TWB}.

Let $\M$ be a manifold of dimension greater than one and let $F:\M \to
\M$ be a mapping. Assume that $\M$ has a product structure so that
$\M=M\times N$ with $M$, $N$ manifolds of dimension at least one, and
$\dim M=m$. By a choice of product structure on $\M$, we mean a choice
of a pair of foliations $\F$ and $\G$ such that each leaf of $\F$ is
homeomorphic to $M$, each leaf of $\G$ is homeomorphic to $N$, and
each leaf of $\F$ crosses each leaf of $\G$ at a single point. If we
choose any leaf $\F _0$ in $\F$, we can now identify $\F _0$ with $M$
and define the map $\Pi _M:\M \to M$ as the projection onto $\F _0$
along $\G$. Similarly, the choice of a leaf $\G _0$ in $\G$ and
an identification of  $\G_0$ with $N$ define a
map $\Pi _N:\M \to N$. We could as well think in terms of choices of
coordinate system and call $\Pi _M (p)$ and $\Pi _N (p)$ the $x$ and
$y$ coordinates of $p\in \M$.  For any choice of product structure, we
can perform the following construction to produce a nonautonomous
mapping on the second factor. Let $\lbrace x(i)\rbrace _{i=0}^\infty$
be a sequence of points in $M$ (in practice, this sequence can be
defined as the successive $x$-coordinates of a sequence of points
$\lbrace p(i)\rbrace _{i=0}^\infty$ in $\M$). Define $\tilde{F}_x:\N
\times N \to N$ by
$$
\tilde{F}_x(n,y)=\Pi_N F(x(n),y) \,.
$$
A nonautonomous mapping on $N$ is defined by
$$
Y(n+1)=\tilde{F}_x(n,Y(n))\,.
$$
With $d_N$ standing for the distance on $N$, we say the system $(\M,
F,\F ,\G )$ is {\it absolutely m-synchronizing} if for all sequences
$\lbrace x(i) \rbrace_{i=0}^\infty$, the motions under $\tilde{F}_x$
on $N$ of all initial conditions converge, {\it i.e.}, for all
$Y_1(0)$, $Y_2(0)$,
$$
\lim_{n\to \infty}d_N(Y_1(n),Y_2(n))=0\,.
$$

If the convergence does not necessarily occur for all sequences
$\lbrace x(i)\rbrace _{i=0}^\infty$, but does occur whenever the
sequence $\lbrace p(i)\rbrace _{i=0}^\infty$ is defined as the
successive points of an orbit of $F$, we say the system $(\M, F,\F ,\G
)$ is {\it m-synchronizing}.  Notice that with $p(i) = (x(i),y(i))$,
for all $Y(0)$ we then have $\lim_{n \rightarrow \infty}d_N(Y(n),y(n))
= 0$, which means that $Y$ synchronizes to $y$.

The map $F$ is {\it (absolutely) m-synchronizable} if a product
structure can be exhibited on $\M$ such that the system $(\M, F,\F ,\G
)$ is (absolutely) m-synchronizing. When $m=1$, the $m$-prefix will
usually be omitted. For instance, any H{\'e}non map is absolutely
synchronizable; in the standard coordinate system, such a map is given
by an equation of the form $F(u,v)=(v,1-av^2+bu)$ so that if one
imposes the time evolution of $v$ then the evolution of $u$ is
independent of $u(0)$ (since $u(n+1) = \tilde{F}_v(n,u(n)) = v(n)$).

One can have good reasons to consider product structures where the
foliations $\F$ and $\G$ have some degree of smoothness, or some other
structure: for instance, particular attention is given in \cite{TWB}
to linear foliations in the case when $\M=\R^d$ and $F$ is a linear
map.
Note that instead of changing the product structure, one can
equivalently keep the product structure fixed and consider all maps
conjugated to $F$, where the smoothness of the conjugacy is the same
as the smoothness of the foliations corresponding to the product
structures one wants to consider.  The following result is
obtained in \cite{TWB}.
\begin{theorem}[\cite{TWB}]
Consider the synchronizing systems $(\R^d,F,\F,\G)$ where $d\geq2$,
$F$ is a linear map or flow, $\F$ is the foliation by lines where
$x_i$ is constant for $i>1$, and $\G$ is the foliation by planes
$x_1=$constant.  Then the set of maps or flows linearly conjugate to
such an $F$ form open and dense sets among linear maps and flows in
$\R ^d$.
\end{theorem}
This theorem admits the following corollary.
\begin{corollary}[\cite{TWB}]
The linear maps and flows in $\R ^d$, $d\geq 2$, which are absolutely
synchronizable form open and dense subsets of the linear maps and
flows on these spaces.
\end{corollary}

\begin{remark}
Consider the system $(\M, F, \F, \G)$.  Suppose the map $F$ has a
fixed point $q_0$.  If we consider the orbit $p(i)\equiv q_0$, then
the map $\tilde{F}_x$ is an autonomous map from the leaf $\F_0$
containing $q_0$ to itself.  Clearly, $q_0$ is a fixed point of this
map.  If $\tilde{F}_x$ has another fixed point, then the system $(\M,
F, \F, \G)$ is not synchronizing.  If this holds true for all pairs of
foliations $\F$, $\G$, then $F$ is not synchronizable.
\end{remark}

Using this remark, it is easy to find examples of diffeomorphisms of
$\R^d$, $d \geq 2$, which are not synchronizable. The identity, and
more generally, any map which leaves the boundary of a topological
$d$-ball pointwise invariant are such maps. However, there are
arbitrarily small perturbations (in any reasonable metric) of these
examples which lose this property. In view of Theorem 1, this leads to
the question of whether non-synchronizability can occur in a robust
way. To make this question more precise, we formulate the following
\begin{definition}
A map is {\it $k$-stably non-synchronizable} if it is $C^k$ and every
map sufficiently $C^k$-close to it is non-synchronizable.  A
homeomorphism is {\it $k$-stably non-synchronizable} if it is $C^k$
and every homeomorphism sufficiently $C^k$ close to it is
non-synchronizable.  When $k=0$, the $k$ prefix will be omitted.
\end{definition}
Examples of linear stably non-synchronizable maps on the
two-dimensional torus $\T ^2$ were given in \cite{TWB} (a more
comprehensive treatment of this matter, including necessary and
sufficient conditions for synchronizability for automorphisms of $\T
^d$, will be given elsewhere \cite{MTW}).  The situation in Euclidean
spaces is quite different.  So far, we could not find any pair
$(r,k)$, $k\leq r$ for which we could exhibit a $C^r$ map in $\R^d$,
$d\geq 3$ which is $k$-stably non-synchronizable. To the contrary in
$\R^2$, we shall prove

\begin{theorem} \label{thm:main}
There exist real-analytic stably non-synchronizable homeomorphisms of $\R^2$.
\end{theorem}

\section{Proof of Theorem \ref{thm:main}}

We denote by $\T ^1=\T$ the one-dimensional torus, and $\pi _1$, $\pi
_2$ the standard projections defined on the plane $\R ^2$ or the
annulus $\A=\T\times [0,1]$ (this ambiguity should not cause any
confusion).  Recall that a simple closed curve $C\subset \A$ is {\it
essential} if it is not null homotopic.  The curve $C'\subset A$ is
{\it above} (respectively, {\it below}) $C$ if it lies in the same
component of $\A \setminus C$ as $\T\times\lbrace 1\rbrace$
(respectively, $\T\times\lbrace 0\rbrace$).  Finally, two points $x$
and $y$ in a topological space $S$ are {\it separated} by a subset
$T\subset S$ if they lie in different connected components of $S
\setminus T$.

\begin{lemma} \label{lemma:1}
Let $F$ be a homeomorphism of $\A$ isotopic to the identity
({\it i.e.}, orientation preserving and leaving each boundary circle
invariant) which satisfies the following conditions:

i) There exists an essential
simple closed curve
$C\subset \T\times (0,1)$ such that $F(C)$ is disjoint from $C$ and
is above $C$, and such that the sequence $(F^n(C))_{n> 0}$
accumulates on $\T\times \lbrace 1\rbrace$ and the sequence
$(F^{-n}(C))_{n > 0}$ accumulates on $\T\times \lbrace 0\rbrace$.

ii) There exists a lift $f$ of $F$ to the universal cover
$\tilde{\A}=\R\times  [0,1] $ of $\A$ such that
$$
\text{for all}\,\, x \in \R\quad \left\{
	\begin{array}{c}
		\pi _1\circ f(x,0)<x\,,\\
		\pi _1\circ f(x,1)>x\,.\\
	\end{array} \right.
$$
Then for each simple arc $\gamma : [0,1] \to\tilde{\A}$ such that
$\gamma (0)\in \R\times \lbrace 0\rbrace$,
$\gamma (1)\in \R\times \lbrace 1\rbrace$,
$\gamma ((0,1))\subset \R\times (0,1)$, there exists $t$ and $t'$
in $(0,1)$ with $t'>t$ such that
$$f(\g (t))=\g (t')\,.$$
\end{lemma}

\begin{proof}[Proof of Lemma \ref{lemma:1}]
We want first to construct a simple arc $\lambda : [0,1] \to \tilde{\A}$
which satisfies the conclusion of the lemma and has the additional
property that there does not exist $u$ and $u'$ in $(0,1)$ with $u'<u$
such that $f(\lambda (u))=\lambda (u')$.

A change of variables may be made in $\R\times (0,1)$ such that each
$f^k(C)$, $k\in \Z$, is of the form $\T\times \lbrace r_k\rbrace$.  It
follows that $\pi ' _2(F(z))>\pi '_2(z)$ for each $z\in \T\times
(0,1)$, where $\pi '_2$ is the second projection with respect to the
new coordinate system.  This fact combined with property $ii)$ implies
that $f$ has no fixed point in $\tilde{\A}$.  Let $C_{\e}$ and
$C_{1-\e}$ be the curves which in the new coordinate system are given
by, respectively, $\R\times \{\e\}$ and $\R\times \{1-\e\}$.  The
curve $C_{\e}$ is pointwise close in the original coordinates to the
boundary component $\R \times \{0\}$ of $\tilde{\A}$, and similarly
for $C_{1-\e}$ to $\R \times \{1\}$.  We now choose one of the points
$p$ of $C_\e$ closest to $\R \times \{0\}$ and $q$ of $C_{1-\e}$
closest to $\R \times \{1\}$.  We define $\lambda$ by sending $[0,\e]$
to the vertical segment in $\tilde{\A}$ joining $\R \times \{0\}$ to
$p$, by sending $[1-\e,1]$ to the vertical segment in $\tilde{\A}$
joining $q$ to $\R \times \{1\}$, and sending $[\e, 1-\e]$ to any arc
joining $p$ to $q$ so that $\pi '_2 \circ \lambda $ is strictly
increasing on $[\e, 1-\e]$.  In the previous sentence, vertical is
with respect to the original coordinate system.  See Figure
\ref{fig:lambda}.
\begin{figure}
\begin{center}
\begin{picture}(184,164)
\put(20,0){\psfig{figure=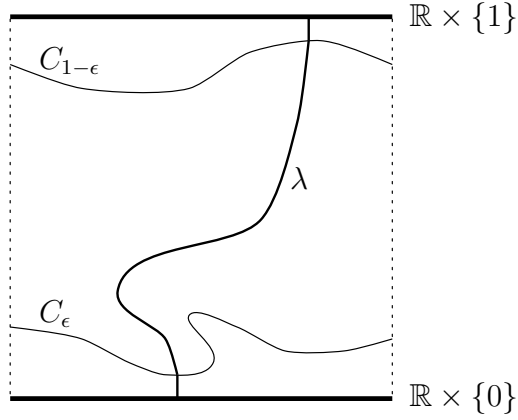}}
\put(170,-5){\makebox(0,0)[lb]{$\R\times\{0\}$}}
\put(170,139){\makebox(0,0)[lb]{$\R\times\{1\}$}}
\put(30,28){\makebox(0,0)[lb]{$C_\e$}}
\put(30,123){\makebox(0,0)[lb]{$C_{1-\e}$}}
\put(125,80){\makebox(0,0)[lb]{$\lambda$}}
\end{picture}
\end{center}
\caption{The curves $C_\e$ and $C_{1-\e}$ and the arc $\lambda$ in the
original coordinate system.} \label{fig:lambda}
\end{figure}
By property $ii)$ there exist $t$ and $t'$ in $(0,1)$ such that
$$
f(\lambda (t))=\lambda (t') \;,
$$
and by the property of $\pi'_2$ necessarily $t'>t$, if
$\e$ is chosen small enough.

Let $\lambda$ be the arc we just constructed and let $\gamma$ be the
arc from the statement of Lemma \ref{lemma:1}.  We can construct a
continuous function
$$
\Gamma : [0,1] \times  [0,1] \to \tilde{\A} \,,
$$
such that, with $\Gamma (s,t)=\g _s(t)$, we have
for each $s$ and $t$ in $ [0,1] $,
$$
\left\{\begin{array}{c}
	{\gamma _0(t)=\lambda(t)\,,} \\
	{\gamma _1 (t)=\gamma (t)\,,}\\
	{\gamma _s(0)\in \R\times \lbrace 0\rbrace\,,}\\
	{\gamma _s(1)\in \R\times \lbrace 1\rbrace\,,}\\
	{\gamma _s (t)\in \R\times (0,1)\,\,\text{if}\,\,t\in (0,1)\,.}\\
\end{array}\right .
$$
We can write $[0,1]^2=U^+\sqcup U^-\sqcup K$ where
$$U^- = \{ (s,t) | f(\gamma _s(t)) \text{ is to the left of }
\gamma _s([0,1])\} \;, $$
$$U^+ = \{ (s,t) | f(\gamma _s(t)) \text{ is to the right of }
\gamma _s([0,1])\}\;, $$
$$ K = \{ (s,t) | f(\gamma _s(t)) \in \gamma _s([0,1]) \} \;.$$
\noindent
Note that $U^-$ and $U^+$ are open subsets of $[0,1]^2$ which
respectively contain $[0,1]\times \lbrace 0\rbrace$ and $[0,1]\times
\lbrace 1\rbrace$, and that $K$ is a closed subset of $[0,1]^2$.

We can define a map $\phi :K\to [0,1]$ by setting
$$f(\gamma _s(t))=\gamma _s (\phi (s,t))\,\,\text{if}\,\,(s,t)\in K\,.$$
This map is continuous and, since $f$ has no fixed point, we can write
$$K=K^+\sqcup K^-\,,$$
where $K^+$ and $K^-$ are closed subsets of $K$ defined by
$$K^+ = \{ (s,t) \in K | \phi(s,t)>t \} \,,$$
$$K^- = \{ (s,t) \in K | \phi(s,t)<t \} \,.$$
Let us remark that
$$
\left\{\begin{array}{c}
	{K^+\cap \lbrace 0 \rbrace\times [0,1]\not =\emptyset\,,}\\
	{K^-\cap \lbrace 0 \rbrace\times [0,1] =\emptyset\,,}
\end{array} \right.
$$
due to the properties of $\lambda$.

To prove Lemma \ref{lemma:1}, we have to show that
$$
K^+\cap \lbrace 1\rbrace \times [0,1]\not =\emptyset\,.
$$
Assume the contrary holds true. Then there exists $\epsilon >0$
such that
$$
\left\{\begin{array}{c}
{K^+\cap [1-\epsilon ,1]\times [0,1] =\emptyset\,,}\\
{K^-\cap [0,\epsilon ]\times [0,1] =\emptyset\,,}
\end{array} \right.
$$
and such that
$$
\left\{\begin{array}{c}
	{[0,1]\times [1-\epsilon ,1]\subset U^+\,,}\\
	{[0,1]\times [0,\epsilon ]\subset U^-\,.}\\ \end{array}
\right.
$$
We set $x=(\frac{1}{2},1-\epsilon )\in U^+$, $y=(\frac{1}{2},
\epsilon )\in U^-$ and we denote by $\Delta$ the boundary of the square
$[0,1]^2$.
Now:
\begin{itemize}
\item the points $x$ and $y$  are separated neither by the closed set
$\Delta \cup K^-$, nor by the closed set $\Delta \cup K^+$,

\item the intersection of $\Delta \cup K^-$ and $\Delta\cup K^+$ is
$\Delta$, which is a connected set.
\end{itemize}
By Alexander's Lemma (see \cite[Theorem 9.2]{Ne}),
we deduce from these two facts that $x$ and $y$ are not separated by
$(\Delta \cup K^-)\cup (\Delta \cup K^+)=\Delta \cup K$.
This is clearly impossible because $x \in U^+$ and $y \in U^-$.
\end{proof}

Returning from working on the universal cover to working on the annulus,
we deduce from Lemma \ref{lemma:1} the following corollary.

\begin{corollary} \label{cor:1}
Under the assumptions of Lemma \ref{lemma:1},
if $\g:[0,1]\to \A$ is an arc starting at $\T \times \lbrace
0\rbrace$, terminating at $\T\times \lbrace 1\rbrace$, and taking its
intermediate values in $\T\times (0,1)$, then there exist $t$ and $t'$
in $(0,1)$ with $t'>t$ such that $F(\g (t))=\g (t')$.
\end{corollary}

If $F$ satisfies the conditions of Lemma \ref{lemma:1},
we shall say that it is {\it of type} (P). We shall
say $F$ is {\it of type} (Q) if it is of type (P) and also satisfies
the following condition $iii)$ which is stronger than condition $ii)$
of Lemma \ref{lemma:1}.

\medskip
{\it
iii) There exists a lift $f$ of $F$
to the universal cover $\tilde \A =\R\times [0,1]$ of $\A$
such that for all $x\in \R$,
$$
\left\{ \begin{array}{c}
	{\pi _1\circ f(x,0)<x-1\,,}\\
	{\pi _1\circ f(x,1)>x+1\,.}
\end{array} \right.
$$
}
\medskip

It is more convenient to present the following fact for the annulus
than for the universal cover.
\begin{lemma} \label{lemma:2}
Assume $F$ is of type (Q).  If $\g$ and $\g '$ are two disjoint arcs
starting at $\T \times \lbrace 0\rbrace$, terminating
at $\T\times \lbrace 1\rbrace$, and taking their intermediate values
in $\T\times (0,1)$, we have
$$
F(\g ([0,1]))\cap \g ' ([0,1])\not = \emptyset \,.
$$
\end{lemma}
\noindent
The proof is straightforward.

The next lemma gives a sufficient condition for non-synchronizability.
Let us consider an orientation preserving homeomorphism $F$ of $\R^2$
which has two compact invariant annuli $\A_1$ and $\A_2$ such that
$\A_1$ is included in the bounded component of $\R ^2\setminus \A_2$.
Assume furthermore that each of $F|_{\A_1}$ and $F|_{\A_2}$ is
conjugated to a homeomorphism of type (Q) by a conjugacy which sends
the inner boundary of $\A_i$ to $\T\times \lbrace 0\rbrace$. We say
such an $F$ satisfies {\it condition} (R).

\begin{lemma} \label{lemma:3}
If $F$ satisfies condition (R), then it is not synchronizable.
\end{lemma}

\begin{proof}[Proof of Lemma \ref{lemma:3}]
By the Brouwer Fixed Point Theorem, $F$ has a fixed point $z_0$ in the
closure of the bounded connected component of $\R ^2\setminus \A_1$.
We can assume that $z_0=(0,0)$. We define
$$\g:[0,+\infty )\to\R\,,\qquad t\mapsto (0,t)\,,$$
and
$$\g ':[0,+\infty )\to\R\,,\qquad t\mapsto (t,0)\,.$$
It is easy to see that there exist two disjoint compact intervals $I_1$
and $I_2$ of $[0,+\infty)$, the first one to the left of the second
one, such that $\g |_{I_i}$ starts at the inner boundary of
$\A_i$, ends at the upper boundary of $\A_i$, and takes its
intermediate values in the interior of $\A_i$. Similarly, there
exists a compact interval $I'_2\subset (0,+\infty)$ such that
$\g '|_{I'_2}$ satisfies the same condition as $\g |_{I_2}$.
See Figure \ref{fig:I}.
\begin{figure}
\begin{center}
\begin{picture}(252,252)
\put(0,0){\psfig{figure=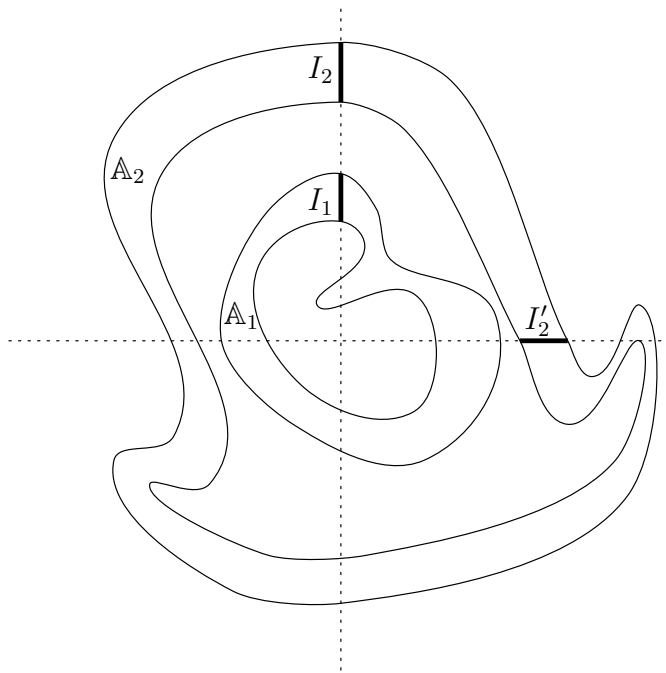}}
\put(45,192){\makebox(0,0){$\A_2$}}
\put(88,137){\makebox(0,0){$\A_1$}}
\put(118,180){\makebox(0,0){$I_1$}}
\put(118,230){\makebox(0,0){$I_2$}}
\put(200,135){\makebox(0,0){$I'_2$}}
\end{picture}
\end{center}
\caption{The annuli $\A_1$ and $\A_2$ and the intervals
$I_1$, $I_2$, and $I'_2$.} \label{fig:I}
\end{figure}
By Corollary \ref{cor:1} there exist $t_1$ and $t_2$ in $I_1$
with $t_1<t_2$, and by Lemma \ref{lemma:2} there exists
$t_3$ in $I_2$ (thus $t_3>t_2$) such that
$$\left\{\begin{array}{c}
{F(\g (t_1))=\g (t_2)\,,}\\
{F(\g (t_3))\in \g '(I'_2)\,.} \end{array} \right . $$
It follows that the continuous map $\psi:t\mapsto \pi _2(F(0,t))$ satisfies
$\psi(t_1)>t_1$ and $\psi (t_3)=0$. Consequently, $\psi$ has a fixed
point between $t_1$ and $t_3$ as well as another fixed point at $0$.
With the notation of the first section, the map $\psi$ is the
same as $\tilde{F}_x$ with $x(i) \equiv 0$.  From the remark in the
first section, $F$ is not synchronizing.
Since condition (R) is stable under conjugacy,
$F$ is not synchronizable.
\end{proof}

\begin{proof}[Proof of Theorem \ref{thm:main}]
Consider now a real analytic homeomorphism $F:\R ^2\to \R ^2$ given in polar
coordinates by $$(\theta, r)\mapsto (\theta+\beta(r), \alpha (r))
\;.$$ We further assume that the map $\alpha :[0,\infty)\to
[0,\infty)$ has fixed points at $r=0,1,2,3,4$ such that 0,2, and 4 are
sinks, while 1 and 3 are sources.

We denote by $S_r$ the disk with center $O$ and radius $r$.  Let
$\A_1$ be the annulus between $S_1$ and $S_2$ and $\A_2$ be the
annulus between $S_3$ and $S_4$.  Finally, we assume that $F|_{\A_1}$
and $F|_{\A_2}$ are of type (Q), so that $F$ satisfies condition (R).
An example of a map $F$ satisfying all these conditions is obtained
by taking
\begin{eqnarray*}
\alpha (r) & = & r+\mu r (r^2-1)(r^2-4)(r^2-9)(r^2-16)(r^2-25) \;, \\
\beta (r) & = & 2\pi r^2 \;,
\end{eqnarray*}
with $\mu$ positive but sufficiently small.

If $G$ is a homeomorphism $C^0$ close to $F$, it does not necessarily
satisfy condition (R) since it might lose the invariant curves. However
we shall see $G$ is not synchronizable, hence $F$ is stably
non-synchronizable. (One could check that
any map close to $F$ in the $C^1$ topology does satisfy condition (R)
and thus get, by Lemma \ref{lemma:3}, an independent proof of the fact that
$F$ is $C^1$-stably non-synchronizable.)

In the rest of the discussion, proximity refers to the $C^0$ norm. If
$G$ is close enough to $F$, the curve $G(S_{3/2})$ is disjoint from
$S_{3/2}$ and situated in the unbounded component of $\R ^2\setminus
S_{3/2}$. For any $n\geq 1$, we can consider the open annulus $V_n$
with boundaries $G^{-n}(S_{3/2})$ and $G^n(S_{3/2})$; these annuli
form an increasing sequence and the set $V=\cup _{n\geq 1}V_n$ is
an open annulus ({\it i.e.}, homeomorphic to $\T\times (0,1)$)
invariant under $G$.

Let us fix $\e\in (0,\frac{1}{2})$. If $G$ is close enough to $F$, the
curves $S_{1-\e}$ and $S_{2+\e}$ are disjoint from their images under
$G$, which are located in the same component, respectively of $\R
^2\setminus S_{1-\e}$ and $\R ^2\setminus S_{2+\e}$ as the images
under $F$. The annulus $V$ is thus bounded between the curves
$G(S_{1-\e})$ and $G(S_{2+\e})$. Furthermore, for $G$ close enough to
$F$, $V$ contains the annulus bounded by $S_{1+\e}$ and $S_{2-\e}$.
Proceeding between $S_3$ and $S_4$ as we just did between $S_1$ and
$S_2$, we get a second open invariant annulus $W$.

By the Brouwer Fixed Point Theorem, we also know that $G$ admits a
fixed point $z_0$ in the closed disk bounded by $S_{1-\e}$. Thus we can
define two maps from $[0,\infty)$ to $\R ^2$:
$$
\g\,:t\mapsto z_0+(0,t)
$$
$$
\g':t\mapsto z_0+(t,0)
$$
and as before, we know that there exist two open intervals $I_1$ and
$I_2$ with disjoint and compact closures, with $I_1$ to the left of
$I_2$ such that $\g |_{I_1}$ and $\g |_{I_2}$ have value respectively
in $V$ and $W$ and join the interior end to the exterior end of these
annuli. Similarly, we know that there exists an open interval $I'_2$
such that $\g '|_{I'_2}$ has the same property as $\g |_{I_1}$.

We now use Caratheodory's prime ends theory (see, {\it e.g.},
\cite[pp. 187-198]{Fo}). We can compactify each connected piece of the
boundary of $V$ as a circle, and extend $G|_V$ as a map $\hat G$ on
the compact annulus $\hat V$ obtained by this compactification. Notice
that $\hat G$ satisfies condition $i)$ with $C=S_{3/2}$. The
interpretation of the prime ends as accessible points tells us that
the rotation numbers of $\hat G$ on the boundary circles of $\hat V$
approach the corresponding rotation numbers of $F|_{S_1}$ and
$F|_{S_2}$ when $G$ converges to $F$; in particular, $\hat G$
satisfies property $iii)$ if $G$ is close enough to $F$.

We have a similar situation for $W$. Furthermore, we know that the
arcs $\g |_{I_1}$, $\g |_{I_2}$, and $\g' |_{I'_2}$ have limits at both
ends which belong to the boundary circles of $\hat V$ or $\hat W$ that
these arcs approach (this is a property of accessible points). Thus we
can apply Corollary \ref{cor:1} and Lemma \ref{lemma:2}, and get the
existence of $t_1<t_2<t_3$ such that
$$
\left\{\begin{array}{c}
	{\hat G(\g (t_1))=\g (t_2)\,,}\\
	{\hat G(\g (t_3))\in \g '(I'_2)\,,}
\end{array} \right.
$$
from which it follows that
$$
\left\{ \begin{array}{c}
{G(\g (t_1))=\g (t_2)\,,} \\
{G(\g (t_3))\in \g '(I'_2)\,.}
\end{array} \right.
$$
We thus conclude that $G$ is not synchronizable.
\end{proof}


\enddocument